\catcode`\@=11

\font\tenmsa=msam10
\font\sevenmsa=msam7
\font\fivemsa=msam5
\font\tenmsb=msbm10
\font\sevenmsb=msbm7
\font\fivemsb=msbm5
\newfam\msafam
\newfam\msbfam
\textfont\msafam=\tenmsa  \scriptfont\msafam=\sevenmsa
  \scriptscriptfont\msafam=\fivemsa
\textfont\msbfam=\tenmsb  \scriptfont\msbfam=\sevenmsb
  \scriptscriptfont\msbfam=\fivemsb

\def\hexnumber@#1{\ifnum#1<10 \number#1\else
 \ifnum#1=10 A\else\ifnum#1=11 B\else\ifnum#1=12 C\else
 \ifnum#1=13 D\else\ifnum#1=14 E\else\ifnum#1=15 F\fi\fi\fi\fi\fi\fi\fi}

\def\msa@{\hexnumber@\msafam}
\def\msb@{\hexnumber@\msbfam}
\mathchardef\boxdot="2\msa@00
\mathchardef\boxplus="2\msa@01
\mathchardef\boxtimes="2\msa@02
\mathchardef\square="0\msa@03
\mathchardef\blacksquare="0\msa@04
\mathchardef\centerdot="2\msa@05
\mathchardef\lozenge="0\msa@06
\mathchardef\blacklozenge="0\msa@07
\mathchardef\circlearrowright="3\msa@08
\mathchardef\circlearrowleft="3\msa@09
\mathchardef\rightleftharpoons="3\msa@0A
\mathchardef\leftrightharpoons="3\msa@0B
\mathchardef\boxminus="2\msa@0C
\mathchardef\Vdash="3\msa@0D
\mathchardef\Vvdash="3\msa@0E
\mathchardef\vDash="3\msa@0F
\mathchardef\twoheadrightarrow="3\msa@10
\mathchardef\twoheadleftarrow="3\msa@11
\mathchardef\leftleftarrows="3\msa@12
\mathchardef\rightrightarrows="3\msa@13
\mathchardef\upuparrows="3\msa@14
\mathchardef\downdownarrows="3\msa@15
\mathchardef\upharpoonright="3\msa@16

\mathchardef\downharpoonright="3\msa@17
\mathchardef\upharpoonleft="3\msa@18
\mathchardef\downharpoonleft="3\msa@19
\mathchardef\rightarrowtail="3\msa@1A
\mathchardef\leftarrowtail="3\msa@1B
\mathchardef\leftrightarrows="3\msa@1C
\mathchardef\rightleftarrows="3\msa@1D
\mathchardef\Lsh="3\msa@1E
\mathchardef\Rsh="3\msa@1F
\mathchardef\rightsquigarrow="3\msa@20
\mathchardef\leftrightsquigarrow="3\msa@21
\mathchardef\looparrowleft="3\msa@22
\mathchardef\looparrowright="3\msa@23
\mathchardef\circeq="3\msa@24
\mathchardef\succsim="3\msa@25
\mathchardef\gtrsim="3\msa@26
\mathchardef\gtrapprox="3\msa@27
\mathchardef\multimap="3\msa@28
\mathchardef\therefore="3\msa@29
\mathchardef\because="3\msa@2A
\mathchardef\doteqdot="3\msa@2B

\mathchardef\triangleq="3\msa@2C
\mathchardef\precsim="3\msa@2D
\mathchardef\lesssim="3\msa@2E
\mathchardef\lessapprox="3\msa@2F
\mathchardef\eqslantless="3\msa@30
\mathchardef\eqslantgtr="3\msa@31
\mathchardef\curlyeqprec="3\msa@32
\mathchardef\curlyeqsucc="3\msa@33
\mathchardef\preccurlyeq="3\msa@34
\mathchardef\leqq="3\msa@35
\mathchardef\leqslant="3\msa@36
\mathchardef\lessgtr="3\msa@37
\mathchardef\backprime="0\msa@38
\mathchardef\risingdotseq="3\msa@3A
\mathchardef\fallingdotseq="3\msa@3B
\mathchardef\succcurlyeq="3\msa@3C
\mathchardef\geqq="3\msa@3D
\mathchardef\geqslant="3\msa@3E
\mathchardef\gtrless="3\msa@3F
\mathchardef\sqsubset="3\msa@40
\mathchardef\sqsupset="3\msa@41
\mathchardef\trianglerighteq="3\msa@44
\mathchardef\trianglelefteq="3\msa@45
\mathchardef\bigstar="0\msa@46
\mathchardef\between="3\msa@47
\mathchardef\blacktriangledown="0\msa@48
\mathchardef\blacktriangleright="3\msa@49
\mathchardef\blacktriangleleft="3\msa@4A
\mathchardef\blacktriangle="0\msa@4E
\mathchardef\triangledown="0\msa@4F
\mathchardef\eqcirc="3\msa@50
\mathchardef\lesseqgtr="3\msa@51
\mathchardef\gtreqless="3\msa@52
\mathchardef\lesseqqgtr="3\msa@53
\mathchardef\gtreqqless="3\msa@54
\mathchardef\Rrightarrow="3\msa@56
\mathchardef\Lleftarrow="3\msa@57
\mathchardef\veebar="2\msa@59
\mathchardef\barwedge="2\msa@5A
\mathchardef\doublebarwedge="2\msa@5B
\mathchardef\angle="0\msa@5C
\mathchardef\measuredangle="0\msa@5D
\mathchardef\sphericalangle="0\msa@5E
\mathchardef\varpropto="3\msa@5F
\mathchardef\smallsmile="3\msa@60
\mathchardef\smallfrown="3\msa@61
\mathchardef\Subset="3\msa@62
\mathchardef\Supset="3\msa@63
\mathchardef\Cup="2\msa@64

\mathchardef\Cap="2\msa@65

\mathchardef\curlywedge="2\msa@66
\mathchardef\curlyvee="2\msa@67
\mathchardef\leftthreetimes="2\msa@68
\mathchardef\rightthreetimes="2\msa@69
\mathchardef\subseteqq="3\msa@6A
\mathchardef\supseteqq="3\msa@6B
\mathchardef\bumpeq="3\msa@6C
\mathchardef\Bumpeq="3\msa@6D
\mathchardef\lll="3\msa@6E

\mathchardef\ggg="3\msa@6F

\mathchardef\circledS="0\msa@73
\mathchardef\pitchfork="3\msa@74
\mathchardef\dotplus="2\msa@75
\mathchardef\backsim="3\msa@76
\mathchardef\backsimeq="3\msa@77
\mathchardef\complement="0\msa@7B
\mathchardef\intercal="2\msa@7C
\mathchardef\circledcirc="2\msa@7D
\mathchardef\circledast="2\msa@7E
\mathchardef\circleddash="2\msa@7F
\def\ulcorner{\delimiter"4\msa@70\msa@70 }
\def\urcorner{\delimiter"5\msa@71\msa@71 }
\def\llcorner{\delimiter"4\msa@78\msa@78 }
\def\lrcorner{\delimiter"5\msa@79\msa@79 }
\def\yen{\mathhexbox\msa@55 }
\def\checkmark{\mathhexbox\msa@58 }
\def\circledR{\mathhexbox\msa@72 }
\def\maltese{\mathhexbox\msa@7A }
\mathchardef\lvertneqq="3\msb@00
\mathchardef\gvertneqq="3\msb@01
\mathchardef\nleq="3\msb@02
\mathchardef\ngeq="3\msb@03
\mathchardef\nless="3\msb@04
\mathchardef\ngtr="3\msb@05
\mathchardef\nprec="3\msb@06
\mathchardef\nsucc="3\msb@07
\mathchardef\lneqq="3\msb@08
\mathchardef\gneqq="3\msb@09
\mathchardef\nleqslant="3\msb@0A
\mathchardef\ngeqslant="3\msb@0B
\mathchardef\lneq="3\msb@0C
\mathchardef\gneq="3\msb@0D
\mathchardef\npreceq="3\msb@0E
\mathchardef\nsucceq="3\msb@0F
\mathchardef\precnsim="3\msb@10
\mathchardef\succnsim="3\msb@11
\mathchardef\lnsim="3\msb@12
\mathchardef\gnsim="3\msb@13
\mathchardef\nleqq="3\msb@14
\mathchardef\ngeqq="3\msb@15
\mathchardef\precneqq="3\msb@16
\mathchardef\succneqq="3\msb@17
\mathchardef\precnapprox="3\msb@18
\mathchardef\succnapprox="3\msb@19
\mathchardef\lnapprox="3\msb@1A
\mathchardef\gnapprox="3\msb@1B
\mathchardef\nsim="3\msb@1C
\mathchardef\napprox="3\msb@1D
\mathchardef\nsubseteqq="3\msb@22
\mathchardef\nsupseteqq="3\msb@23
\mathchardef\subsetneqq="3\msb@24
\mathchardef\supsetneqq="3\msb@25
\mathchardef\subsetneq="3\msb@28
\mathchardef\supsetneq="3\msb@29
\mathchardef\nsubseteq="3\msb@2A
\mathchardef\nsupseteq="3\msb@2B
\mathchardef\nparallel="3\msb@2C
\mathchardef\nmid="3\msb@2D
\mathchardef\nshortmid="3\msb@2E
\mathchardef\nshortparallel="3\msb@2F
\mathchardef\nvdash="3\msb@30
\mathchardef\nVdash="3\msb@31
\mathchardef\nvDash="3\msb@32
\mathchardef\nVDash="3\msb@33
\mathchardef\ntrianglerighteq="3\msb@34
\mathchardef\ntrianglelefteq="3\msb@35
\mathchardef\ntriangleleft="3\msb@36
\mathchardef\ntriangleright="3\msb@37
\mathchardef\nleftarrow="3\msb@38
\mathchardef\nrightarrow="3\msb@39
\mathchardef\nLeftarrow="3\msb@3A
\mathchardef\nRightarrow="3\msb@3B
\mathchardef\nLeftrightarrow="3\msb@3C
\mathchardef\nleftrightarrow="3\msb@3D
\mathchardef\divideontimes="2\msb@3E
\mathchardef\varnothing="0\msb@3F
\mathchardef\nexists="0\msb@40
\mathchardef\mho="0\msb@66
\mathchardef\thorn="0\msb@67
\mathchardef\beth="0\msb@69
\mathchardef\gimel="0\msb@6A
\mathchardef\daleth="0\msb@6B
\mathchardef\lessdot="3\msb@6C
\mathchardef\gtrdot="3\msb@6D
\mathchardef\ltimes="2\msb@6E
\mathchardef\rtimes="2\msb@6F
\mathchardef\shortmid="3\msb@70
\mathchardef\shortparallel="3\msb@71
\mathchardef\smallsetminus="2\msb@72
\mathchardef\thicksim="3\msb@73
\mathchardef\thickapprox="3\msb@74
\mathchardef\approxeq="3\msb@75
\mathchardef\succapprox="3\msb@76
\mathchardef\precapprox="3\msb@77
\mathchardef\curvearrowleft="3\msb@78
\mathchardef\curvearrowright="3\msb@79
\mathchardef\digamma="0\msb@7A
\mathchardef\varkappa="0\msb@7B
\mathchardef\hslash="0\msb@7D
\mathchardef\hbar="0\msb@7E
\mathchardef\backepsilon="3\msb@7F
\def\Bbb{\ifmmode\let\next\Bbb@\else
 \def\next{\errmessage{Use \string\Bbb\space only in math mode}}\fi\next}
\def\Bbb@#1{{\Bbb@@{#1}}}
\def\Bbb@@#1{\fam\msbfam#1}

\catcode`\@=\active

\def\inv{^{\raise.15ex\hbox{${
  \scriptscriptstyle -}$}\kern-.05em 1}}

\def\Dsl{\,\raise.15ex\hbox{$/$}\mkern-13.5mu D}
\def\dsl{\raise.15ex\hbox{$/$}\kern-.57em\hbox{$\partial$}}

\def\lspace{\ifx\answ\bigans{}\else\qquad\fi}

\def\cg{\hbox{{\sl g}}} % used for Lie algebra 'gothic g'

\def\lform{\hbox{$\sqcup$}\llap{\hbox{$\sqcap$}}}
\def\darr#1{\raise1.5ex\hbox{$\leftrightarrow$}
\mkern-16.5mu #1}
\def\INT{{\textstyle \int\kern-.642em\int}}

\def\C{{\Bbb C}}
\def\Z{{\Bbb Z}}

\def\eps{{\epsilon}}

\def\lcross{{>\!\!\!\triangleleft}}
\def\rcocross{{\blacktriangleright\!\!<}}

\def\bicross{{\blacktriangleright\!\!\!\triangleleft}}

\def\rbiprod{{\cdot\kern-.33em\triangleright\!\!\!<}}
\def\lbiprod{{>\!\!\!\triangleleft\kern-.33em\cdot}}

\def\tens{\mathop{\otimes}}

\def\la{{\triangleright}}\def\ra{{\triangleleft}}

\def\isom{{\cong}}

\def\id{{\rm id}}

\def\<{\langle}
\def\>{\rangle}

\def\vecl{{\bf l}}

\def\vecm{{\bf m}}

\def\<{\langle}
\def\>{\rangle}

\def\equad{\kern -1.7em}

\def\o{{}_{\scriptscriptstyle(1)}}
\def\t{{}_{\scriptscriptstyle(2)}}
\def\th{{}_{\scriptscriptstyle(3)}}

\def\bo{{}^{\bar{\scriptscriptstyle(1)}}}
\def\bt{{}^{\bar{\scriptscriptstyle(2)}}}

\def\text#1{\mbox{\rm #1}}
\def\note#1{}

\def\blacksquare{{\lform}}%AMS Tex Fakes
\def\frac#1#2{{{#1\over#2}}}

\def\proof{\goodbreak\noindent{\bf Proof\quad}}

\def\endproof{{\ $\lform$}\bigskip }

\def\eqn#1#2{\begin{equation}#2\label{#1}\end{equation}}
\def\align#1{\begin{eqnarray*}#1\end{eqnarray*}}
\def\ceqn#1#2{\begin{equation}\label{#1}\begin{array}{c}#2\end{array}
\end{equation}}

\def\vecM{{\bf M}}

\documentstyle[11pt]{article}
\newtheorem{lemma}{Lemma}[section] \newtheorem{propos}[lemma]{Proposition}

\begin{document}\baselineskip 22pt

{\ }\qquad\qquad \hskip 4.3in DAMTP/95-59
\vspace{.2in}

\begin{center} {\LARGE BICROSSPRODUCT STRUCTURE OF AFFINE\\ QUANTUM GROUPS}
\\ \baselineskip 13pt{\ }
{\ }\\ S. Majid\footnote{Royal Society University Research Fellow and Fellow of
Pembroke College, Cambridge. On leave 1995 + 1996 at the Department of
Mathematics, Harvard University, Cambridge MA02138, USA}\\
{\ }\\
Department of Applied Mathematics \& Theoretical Physics\\
University of Cambridge, Cambridge CB3 9EW\\
+\\
Research Institute of Mathematical Sciences\\
Kyoto University, Kyoto 606, Japan
\end{center}
\begin{center}
July -- revised November, 1995
\end{center}

\vspace{10pt}
\begin{quote}\baselineskip 13pt
\noindent{\bf Abstract}
We show that the affine quantum group $U_q(\hat{sl_2})$ is isomorphic to a
bicrossproduct central extension $\C\Z{}_\chi\bicross U_q(Lsl_2)$ of the
quantum loop group $U_q(Lsl_2)$ by a quantum cocycle $\chi$, which we
construct.  We prove the same result for $U_q(\hat{\cg})$ in R-matrix form.


\end{quote}
\baselineskip 22pt

\section{Introduction}

Affine quantum groups figure prominently as the quantum `non-Abelian
symmetries' of q-deformed conformal field theories and certain statistical
models\cite{JimMiw:alg}. Their representation theory has been extensively
studied via the techniques of vertex
algebras\cite{FreRes:aff}. A result which has been missing, however, is the
sense in which these quantum groups are central extensions of quantum loop
groups. We recall that  this is important for the correct geometrical picture
in the classical theory \cite{PreSeg}, hence should be important for a
geometrical picture in the $q$-deformed case as well.

We provide this result in the present paper, constructing the appropriate
`quantum cocycle' $\chi:U_q(L\cg)^{\tens 2}\to \C\Z$, from which
$U_q(\hat{\cg})$ is then obtained as the corresponding extension. Here
$U_q(L\cg)$ is the level 0 version of the affine quantum group, and $\C\Z$
denotes the group algebra of $\Z$, i.e. polynomials in a generator $c$ and its
inverse.

The required theory of (non-Abelian) Hopf algebra extensions was already
introduced (by the author) in \cite{Ma:mor}, where it is shown that extensions
generally have the form of a {\em cocycle bicrossproduct}.  The quantum Weyl
group is already known to be of this form\cite{MaSoi:bic}. In general however,
it can be hard to come up with the cocycle data required in the
construction. Our first result, in Section~2, provides a general solution to
the contruction problem in the case of a central extension.

We give the case of $U_q(\hat{sl_2})$ in detail, in Section~3. In Section~4 we
give the result more generally using the  $R$-matrix formalism with generators
$\vecl^\pm(z)$ in \cite{ResSem:cen}\cite{FreRes:aff}.

\subsection*{Acknowledgements} These results were obtained during a visit in
June 1995 to R.I.M.S. in Kyoto under a joint programme with the Isaac Newton
Institute in Cambridge and the J.S.P.S. I would like to thank M. Jimbo for
some useful discussions.

\section{General Construction of Quantum Cocycles}

In this section we introduce a general construction for quantum cocycles
appropriate to the central extensions which concern us. We begin by recalling
the more abstract (but not constructive) theory of cocycle bicrossproducts and
Hopf algebras from \cite{Ma:mor}\cite{MaSoi:bic} in the form now required.

We assume standard notations for Hopf algebras $H,A$, coproducts $\Delta
h=h\o\tens h\t$, antipodes $S$, left and
right actions $\la,\ra$, coaction $\beta(h)=h\bo\tens h\bt$, etc.; see the
texts\cite{Swe:hop}\cite{Ma:book}.
 A module algebra means an algebra
which is acted upon covariantly. Likewise, a comodule coalgebra means a
coalgebra which is coacted
upon covariantly.

Let $H$ be a Hopf algebra and $A$ an algebra. A quantum 2-cocycle on $H$ with
values in $A$ means a map $\chi:H\tens H\to A$ obeying
\eqn{cocy}{\chi(g\o \tens f\o)\chi(h\tens g\t  f\t)=\chi(h\o\tens g\o)\chi(h\t
g\t\tens f),\quad
\chi(1\tens (\ ))=\chi((\ )\tens 1)=\eps.}
It is also possible to allow here an action of $H$ on $A$, but we will not need
this. Such cocycles are known to allow\cite{Doi:equ} the formation of an
algebra structure $A{}_\chi\lcross H$ with product
\eqn{extprod}{ (a \tens h)(b\tens g)=a b \chi(h\o\tens g\o)\tens h\t g\t.}
To this older theory, we add\cite{Ma:mor} the supposition that $\beta:H\to
H\tens A$ makes $H$ a right comodule coalgebra compatible with $\chi$ in the
sense
\eqn{ext1}{\beta(hg)=(1\tens \chi^{-1}(h\o\tens
g\o))\beta(h\t)\beta(g\t)(1\tens \chi(h\th\tens g\th)),\quad \beta(1)=1\tens 1}
\eqn{ext2}{[\beta(h),1\tens a]=0}
\eqn{ext3}{ \Delta \chi(h\tens g)=\chi(h\o\bo\tens g\o\bo)\tens h\o\bt g\o\bt
\chi(h\t\tens g\t),\quad \eps\chi
=\eps\tens\eps.}
In this case the standard cross coproduct coalgebra $A\rcocross H$, with
coproduct
\eqn{extcoprod}{ \Delta (a\tens h)=a\o\tens h\o\bo\tens a\t h\o\bt\tens h\t}
forms with the algebra (\ref{extprod}) a Hopf algebra, denoted
$A{}_\chi\bicross H$, the cocycle bicrossproduct. Projection to $H$ by applying
the counit $\eps$ of $A$, and the inclusion of $A$ as $A\tens 1$ provide Hopf
algebra maps
\eqn{exthopf}{ A\hookrightarrow A{}_\chi\bicross H\to H}
obeying certain properties. It is possible to characterise cocycle
bicrossproducts more abstractly as Hopf algebras  $A\hookrightarrow  E\to H$
with the maps obeying certain properties making $E$ an abstract (cleft and
cocleft) extension of Hopf algebras. Then $E\isom A{}_\chi\bicross^\psi H$ of
the general type in \cite{Ma:mor} with, possibly, an action of $H$ on $A$ and a
dual cocycle $\psi$ in the coalgebra. The special case
(\ref{cocy})--(\ref{extcoprod}) has the additional features that $A$ commutes
with the elements of $H$ and application of the counit of $H$ is a coalgebra
map $E\to A$.

There remains the problem of how to actually construct $\chi,\beta$ obeying
(\ref{cocy}), (\ref{ext1})--(\ref{ext3}). Part of this, the coaction $\beta$,
is  usually easy to identify. The following lemma shows that once this is done,
the remaining cocycle can be found, and provides a formula for it.

\begin{propos} Let $A,E,H$ be Hopf algebras with $E=A\rcocross H$ as a
coalgebra
for a coaction $H\to H\tens A$. Suppose that $E\to H$ by the counit of $H$ and
$A\subset E$ by the canonical inlcusion are algebra homomorphisms and that the
linear map $j:H\to E$ defined by $j(h)=1\tens h$ obeys $aj(h)=a\tens h=j(h)a$
for all $h\in H$ and $a\in A$. Then
\[ \chi(h\tens g)=j(h\o)j(g\o)(Sj((h\t g\t)\bo)) (h\t g\t)\bt\]
is a quantum 2-cocycle on $H$ with values in $A$, and $E=A{}_\chi\bicross H$.
\end{propos}
\proof We assume that $H$ is a right $A$-comodule coalgebra so that we can form
the coproduct (\ref{extcoprod}) on $A\tens H$. Then both $A\subset E$ and
$\pi:E\to H$ defined by $a\mapsto a\tens 1$ and $a\tens h\mapsto \eps(a)h$
become Hopf algebra maps. There is also a coaction $\Delta_R:E\to E\tens H$
given by $(\id\tens\pi)\circ\Delta$, such that $A=E^H=\{e\in
E|\Delta_R(e)=e\tens 1\}$, the fixed point subalgebra. We show now that the
linear map $j$ is convolution-invertible, i.e. we find a linear map
$j^{-1}:H\to E$ such that $j(h\o)j^{-1}(h\t)=\eps(h)=j^{-1}(h\o)j(h\t)$ for all
$h\in H$. Indeed, we set
\[ j^{-1}(h)=(S j(h\bo))h\bt\]
where $S$ denotes the antipode of $E$ and the products are in $E$. Then
\[j^{-1}(h\o)j(h\t)=(Sj(h\o\bo))h\o\bt j(h\t)=(S(1\tens h\o\bo))(h\o\bt\tens
h\t)=\eps(1\tens h)=\eps(h)\]
using the form of the coproduct (\ref{extcoprod}) of $E$. On the other side, we
have
\align{&&\equad j(h\o)j^{-1}(h\t)=(1\tens h\o)(S(1\tens h\t\bo) )h\t\bt\\
&&=(1\tens h\o\bo)(S(1\tens h\t\bo))(Sh\o\bt\o)h\o\bt\t h\t\bt\\
&&=(1\tens h\o\bo\bo)(S(1\tens h\t\bo))(Sh\o\bo\bt)h\o\bt h\t\bt\\
&&=(1\tens h\bo\o\bo)(S(1\tens h\bo\t))(Sh\bo\o\bt)h\bt\\
&&=(1\tens h\bo\o\bo)(S(h\bo\o\bt\tens h\bo\t)) h\bt=\eps(h\bo)h\bt=\eps(h)}
where the second equality uses the antipode axiom in $A$ to insert
$(Sh\o\bt\o)h\o\bt\t$, the third  uses the coaction axiom, the fourth   uses
covariance of the coproduct of $H$ under the coaction and the fifth that $A$ is
a sub-Hopf algebra. We can then use that $S$ is the antipode in $E$ to collapse
the expression.
Next, from the form of the coproduct it is clear that the map $j$ intertwines
the above coaction $\Delta_R$ and the right regular coaction of $H$ on itself.
Hence all the conditions for a cleft extension are satisfied and we know from
general theory of algebra extensions\cite{BlaMon:cro} that $\chi(h\tens
g)=j(h\o)j(g\o)j^{-1}(h\t g\t)$ has values in $A$ and obeys (\ref{cocy}) if
$aj(h)=j(h)a$. This is also easy enough to verify directly along the lines
above. We can now  make the cocycle product (\ref{extprod}) and identify this
as the algebra of $E$. We do not need to verify   (\ref{ext1})--(\ref{ext3})
directly since we already know that $E$ is a Hopf algebra; these conditions are
equivalent to the bialgebra homomorphism property \cite{Ma:mor}. \endproof

Extensions of this type can be viewed as quantum principal bundles in the sense
of \cite{BrzMa:gau}, see \cite{Ma:non} for a discussion. We do not require that
$A$ is actually commutative, but when this is so it appears in the center of
$E$ in the setting  above.

\section{$U_q(\hat{sl_2})$ as a quantum central extension}

In this section we show that the affine quantum group $U_q(\hat{sl_2})$ in
\cite{Dri} is a quantum group central extension. This quantum group plays a
central role in
an approach to certain quantum statistical systems\cite{JimMiw:alg}. In the
conventions of the latter, we have generators $E_i,K_i,F_i$,where $i=0,1$ mod
$2$, and relations\cite{Dri}
\ceqn{usl2aff}{ K_iE_i=qE_iK_i,\quad K_iF_i=q^{-1}F_iK_i,
\quad K_iE_{i+1}=q^{-1}E_{i+1}K_i, \quad K_iF_{i+1}=qF_{i+1}K_i\\K_iK_j=K_jK_i,
\quad [E_i,F_j]=\delta_{ij}{K_i-K_i^{-1}\over q-q^{-1}},\\
{}[E_0^3,E_1]=[3]_q E_0[E_0,E_1]E_0,
\quad [F_0^3,F_1]=[3]_q F_0[F_0,F_1]F_0,}
where $[n]_q={q^n-q^{-n}\over q-q^{-1}}$, and the coproduct
\eqn{usl2cop}{\Delta K_i=K_i\tens K_i,\quad \Delta E_i=E_i\tens K_i+1\tens E_i,
\quad \Delta F_i=F_i\tens 1+K_i^{-1}\tens F_i.}

Our first step is to work with new generators $K_iF_i\mapsto F_i$, $K_1\to K$
and $K_0\mapsto cK^{-1}$ so that the relations become $c$ central and
\ceqn{newsl2}{ K E_0=q^{-1}E_0K,\quad KE_1=qE_1K,\quad KF_0=qF_0K,\quad
KF_1=q^{-1}F_1K,\\
 qE_0F_0-F_0E_0={K^{-2}c^2-1\over q-q^{-1}},\quad qE_1F_1-F_1E_1={K^2-1\over
q-q^{-1}}\\
\Delta K=K\tens K,\quad \Delta E_0=E_0\tens cK^{-1}+1\tens E_0,\quad  \Delta
E_1\tens K+1\tens E_1\\
\Delta c=c\tens c,\quad
\Delta F_0=F_0\tens cK^{-1}+1\tens F_0,\quad  \Delta F_1\tens K+1\tens F_1,}
along with the $q$-Serre relations of the same form as in (\ref{usl2aff}) in
terms of our new generators $E_i,F_i$; the $c,K$ cancel from these.

Let $U_q(Lsl_2)$ denote the quotient of $U_q(\hat{sl_2})$ obtained by setting
$c=1$. We let $e_i,f_i,k$ denote the images of the generators $E_i,F_i,K$.
Thus, $qe_0f_0-f_0e_0={k^{-2}-1\over q-q^{-1}}$ and $\Delta e_0=e_0\tens
k^{-1}+1\tens e_0$, etc. Since $c$ is central and group-like, this quotient
remains a Hopf algebra. It is clear that we have Hopf algebra maps
\eqn{projsl2}{ \C\Z\hookrightarrow U_q(\hat{sl_2})\to U_q(Lsl_2)}
where $\C\Z=\C[c,c^{-1}]$ is the sub-Hopf algebra generated by $c,c^{-1}$. By
moving the $c$'s to the left, we can clearly identify
$U_q(\hat{sl_2})=\C\Z\tens U_q(Lsl_2)$ as linear spaces. This identification
restricted to $U_q(Lsl_2)$ is the linear map
\eqn{jsl2}{ j(k^n g(e_0,e_1)h(f_0,f_1))=K^n g(E_0,E_1)h(F_0,F_1),\quad
j:U_q(Lsl_2)\to U_q(\hat{sl_2}).}
Here a general element of $U_q(Lsl_2)$ is clearly a linear combination of terms
of this type (i.e the quantum group has a triangular decomposition), where $g$
is a polynomial in the non-commuting generators $e_i$, and $h$ in the $f_i$.

Next, because the $q$-Serre relations in (\ref{usl2aff}) are homogeneous in the
$E_0,F_0$, the same is true for the $e_0,f_0$. Hence we have a well-defined
$\Z$-grading defined for polynomials $g,h$ which are homogeneous in $e_0,f_0$
respectively. We write $|g|,|h|$ for the total degree of $e_0,f_0$ respectively
(i.e. each term in $g$ has $|g|$ occurences of $e_0$, etc.). From this grading,
we can define a coaction
\eqn{betasl2}{ \beta(k^n g(e_0,e_1)h(f_0,f_1))=k^n g(e_0,e_1)h(f_0,f_1)\tens
c^{|g|+|h|} ,\quad \beta:U_q(Lsl_2)\to U_q(Lsl_2)\tens \C\Z}
on homogeneous polynomials. This is a coaction because applying it again gives
$c^{|g|+|h|}\tens c^{|g|+|h|}=\Delta c^{|g|+|h|}$. Note that it does not
respect the algebra structure, so (as for the map $j$) it is not enough to give
it on generators. It does, however, respect the coproduct because the coproduct
preserves our above ordering (in which we write any $e_i$ to the left of any
$f_i$), and manifestly preserves the $e_0$ and $f_0$ degrees when acting on
generators. Hence the product of the coactions on $\Delta (k^n
g(e_0,e_1)h(f_0,f_1))$ gives the same result as
applying the coaction first as in (\ref{betasl2}) and then $\Delta$.

We are therefore in a position to make the cross coproduct coalgebra as in
(\ref{extcoprod}). For example,
\[ \Delta (1\tens e_0)=(1\tens e_0)\tens (c\tens k^{-1})+ (1\tens 1)\tens
(1\tens e_0)\]
from (\ref{extcoprod}). This is the coproduct of $U_q(\hat{sl_2})$ on its
identification with $\C\Z\tens U_q(Lsl_2)$. Similarly for $\Delta (1\tens
f_0)$. Hence we are in the setting of Proposition~2.1. We conclude:

\begin{propos} $U_q(\hat{sl_2})$ has the structure of a cocycle bicrossproduct
$\C\Z{}_\chi\bicross U_q(Lsl_2)$ with cross coproduct using coaction $\beta$
from (\ref{betasl2}) and the cocycle product
(\ref{extprod}) with respect to
$\chi:U_q(Lsl_2)^{\tens 2}\to \C\Z$ defined by Proposition~2.1 and $j$ from
(\ref{jsl2}).
\end{propos}

For example,
\align{&&\equad \chi(e_0^a\tens f_0^b)=\sum_{r=0,s=0}^{r=a,s=b} \left[{a\atop
r};q^{-1}\right]\left[{b\atop s};q\right] j(e_0^r)j(f_0^s)
(Sj(k^{-r}e_0^{a-r}k^{-s}f_0^{b-s}))c^{a+b-r-s}\\
&&=\sum_{r=0,s=0}^{r=a,s=b} \left[{a\atop r};q^{-1}\right]\left[{b\atop
s};q\right] E_0^rF_0^s
(S(K^{-r}E_0^{a-r}K^{-s}F_0^{b-s}))c^{a+b-r-s}=\delta_{a,0}\delta_{b,0}}
since we arrive at $\cdot\circ(\id\tens S)\circ\Delta(E_0^aF_0^b)$ in this
case. Here $\left[{a\atop r};q^{-1}\right]$ denotes the appropriate
$q$-binomial coefficient in $\Delta e_0^a$, etc. However,
\align{&&\equad \chi(f_0\tens
e_0)=F_0E_0K^2+F_0(Sj(k^{-1}e_0))c+E_0(Sj(f_0k^{-1}))c+(Sj((f_0e_0)\bo))
(f_0j_0)\bt\\
&&=F_0E_0K^2+F_0 cK SE_0  + E_0cK SF_0  + (Sj(qe_0f_0 ))c^2+Sj({1-K^{-2}\over
q-q^{-1}})\\
&&=F_0E_0K^2+F_0 cK SE_0  + E_0cK SF_0  +  c^2 S(qE_0F_0 )+{1-K^2\over
q-q^{-1}}\\
&&=F_0E_0K^2+F_0 cK SE_0  + E_0cK SF_0  +  c^2 S(F_0E_0 ) +{K^2-c^2\over
q-q^{-1}}+{1-K^2\over q-q^{-1}}\\
&&={1-c^2\over q-q^{-1}}.}
Similarly for more general $\chi(f_0^a\tens e_0^b)$. This is the nontrivial
part of the cocycle $\chi$.

One can check that the cocycle product from (\ref{extprod}) indeed recovers the
correct one for $U_q(\hat{sl_2})$. For example,
\[ F_0E_0=\chi(f_0\tens e_0)K^{-2}+ j(f_0e_0)={1-c^2\over
q-q^{-1}}K^{-2}+j(qe_0f_0)+{1-K^{-2}\over q-q^{-1}}=qE_0F_0+{1-c^2K^{-2}\over
q-q^{-1}}\]
as required.

\section{$R$-matrix form of the quantum cocycle}

Here we show that $U_q(\hat{\cg})$ is likewise a quantum group central
extension, at least for those cases which can be treated using the $R$-matrix
formalism in \cite{ResSem:cen}\cite{FreRes:aff}. These authors
identified suitable generators $\vecl^\pm(z),c$, and the relations
\ceqn{ugaff}{ \vecl^\pm_1(z)\vecl_2^\pm(w)R({z\over w})
=R({z\over w})\vecl^\pm_2(w)\vecl^\pm_1(z),\quad
\vecl^-_1(z)\vecl_2^+(w)R({z\over w}q^c)
=R({z\over w}q^{-c})\vecl^+(w)_2\vecl^-_1(z)\\
{}[c,\vecl^\pm(z)]=0,\quad \Delta q^c=q^c\tens q^c,\quad \Delta \vecl^\pm(z)
=\vecl^\pm(zq^{\pm c_2/2})\tens \vecl^\pm(zq^{\mp c_1/2})}
where $c_1=c\tens 1$ and $c_2=1\tens c$ and $R(z)$ is a suitable solution of
the parametrized quantum Yang-Baxter equation $R_{12}({z\over
w})R_{13}(z)R_{23}(w)=R_{23}(w)R_{13}(z)R_{12}({z\over w})$. There are still
further relations obeyed by the $\vecl^\pm$ which, together with
the above `quadratic' ones, provide a definition of $U_q(\hat{\cg})$ in these
terms. Note that the conventions we
use here are not quite those in \cite{ResSem:cen} but are more in line with the
established  conventions
for $U_q(\cg)$ in \cite{FRT:lie}.

First, we move to new matrix generators
\eqn{mpm}{ \vecM^\pm(z)=\vecl^\pm(zq^{\pm {c\over 2}})}
so that the mixed relations and coproduct become
\eqn{mixedmpm}{ \vecM_1^-(z)\vecM_2^+(w)R({z\over w})=R({z\over
w}q^{-2c})\vecM^+_2(w)\vecM^-_1(z),\quad \Delta\vecM^\pm(z)=\vecM^\pm(zq^{\pm
c_2})\tens \vecM^\pm(z).}
The $\vecM^\pm_1\vecM^\pm_2$ relations are unaffected. The antipode is
\eqn{antmpm}{S\vecM^\pm(z)=(\vecM^\pm)^{-1}(zq^{\mp c})}
where $(\vecM^\pm)^{-1}(z)$ is the `pointwise' inverse  matrix-valued
powerseries to $\vecM^\pm(z)$.

Let $U_q(L\cg)$ be the quotient of $U_q(\hat{\cg})$ obtained by setting
$q^c=1$. This is a Hopf
algebra since $q^c$ is grouplike and central. We denote the matrix generators
in this quotient by $\vecm^\pm$. We have Hopf algebra maps
\eqn{extg}{ \C\Z\hookrightarrow U_q(\hat{\cg})\to U_q(L\cg).}
Here $\C\Z$ is the group algebra of $\Z$ with group-like
generator $q^c$ (or the Lie algebra $U(1)$ with primitive $c$). It is clear
that we can identify $U_q(\hat{\cg})=\C\Z\tens U_q(L\cg)$ as linear spaces by
putting $q^c$ to the left in all normal-ordered expressions. Normal ordered
means for us that all $\vecM^-$ modes are put to the right of all $\vecM^+$
modes using the cross relations (\ref{mixedmpm}). Once ordered correctly, we
identify products of $\vecM^\pm$ are corresponding to products of $\vecm^\pm$
in $U_q(L\cg)$, i.e. restriction of this identification gives the linear map
$j:U_q(L\cg)\to U_q(\hat{\cg})$,
\eqn{jmpm}{
j(\vecm^+_1(z_1)\cdots\vecm^+_i(z_i)\vecm^-_{i+1}(z_{i+1})\cdots\vecm^-_j(z_j))
=\vecM^+_1(z_1)\cdots\vecM^+_i(z_i)\vecM^-_{i+1}(z_{i+1})\cdots\vecM^-_j(z_j).}

Similarly, we define a coaction $U_q(L\cg)\to U_q(L\cg)\tens\C\Z$,
\eqn{betampm}{
\beta(\vecm^+_1(z_1)\cdots\vecm^+_i(z_i)\vecm^-_{i+1}(z_{i+1})\cdots
\vecm^-_j(z_j))=\vecm^+_1(z_1q^{c_2})\cdots\vecm^+_i(z_iq^{c_2})\vecm^-_{i+1}
(z_{i+1}q^{-c_2})\cdots\vecm^-_j(z_jq^{-c_2}),}
where $c_2=1\tens c$. It is easy to see that this is a coaction and respects
the matrix form $\Delta \vecm^\pm(z)=\vecm^\pm(z)\tens\vecm^\pm(z)$ of the
coproduct of $U_q(L\cg)$ (making it a comodule coalgebra).

The cross coproduct coalgebra structure by this coaction recovers the coproduct
of $U_q(\hat{\cg})$. Thus,
\[ \Delta \vecM^\pm(z)=\Delta(1\tens\vecm^\pm(z))=(1\tens \vecm^\pm(zq^{\pm
c_3}))\tens (1\tens \vecm^\pm(z))=\vecM^\pm(zq^{\pm c_2})\tens\vecM^\pm(z),\]
where $c_3=1\tens 1\tens c\tens 1$ in the middle expression. Hence we are in
the setting of Proposition~2.1 and conclude:

\begin{propos} $U_q(\hat{\cg})$ has the structure of a cocycle bicrossproduct
$\C\Z{}_\chi\bicross U_q(L\cg)$
with cross coproduct from (\ref{extcoprod}) using coaction (\ref{betampm}) and
cocycle product (\ref{extprod}) with respect to $\chi:U_q(L\cg)^{\tens 2}\to
\C\Z$ defined from Proposition~2.1 and $j$ from (\ref{jmpm}).
\end{propos}

For example, we have
\[\chi(\vecm^+_1(z)\tens \vecm^-_2(w))=
\vecM^+_1(z)\vecM^-_2(w)S(\vecM^+_1(zq^{c_2})\vecM^-_2(wq^{-c_2}))=\id,\]
where $c_2$ denotes $c$ placed to the far right (outside the range of $S$).
These factors cancel $q^{\pm c}$ in (\ref{antmpm}) when we compute $S$, giving
the identity matrix. In the other order, we have
\align{&&\equad \chi(\vecm^-_1(z)\tens \vecm^+_2(w)) =
\vecM^-_1(z)\vecM^+_2(w)(Sj((\vecm^-_1(z)\vecm^+_2(w))\bo))
(\vecm^-_1(z)\vecm^+_2(w))\bt\\
&&=\vecM^-_1(z)\vecM^+_2(w)Sj(R({z\over
w})\vecm^+_2(wq^{c_2})\vecm^-_1(zq^{-c_2})R^{-1}({z\over w}))\\
&&=\vecM^-_1(z)\vecM^+_2(w)R({z\over
w})(\vecM^-)_2^{-1}(z)(\vecM^+)_1^{-1}(w)R^{-1}({z\over w})=R({z\over
w}q^{-2c})R^{-1}({z\over w})}
where we used the relations (\ref{mixedmpm}) to normal order before applying
the coaction $\beta$ and the map $j$. The effect of the coaction cancels the
$q^{\mp c}$ factor from the action of $S$. The final step uses the relations
(\ref{mixedmpm}) in reverse.

It should be clear that the same computation works in general and gives us
$\chi$ on a general product of $\vecm^\pm$ in terms  of products of $R$ and
$R^{-1}$ with $q^{\pm c}$ in the arguments. For example,
\align{&&\equad
\chi(\vecm^-_1(z_1)\vecm^-_2(z_2)\tens\vecm^+_3(z_3)\vecm^+_4(z_4))\\ &&
=R_{23}({z_2\over z_3}q^{-2c})R_{13}({z_1\over z_3}q^{-2c})R_{24}({z_2\over
z_4}q^{-2c})R_{14}({z_1\over z_4}q^{-2c})R_{14}^{-1}({z_1\over
z_4})R_{24}^{-1}({z_2\over z_4})R_{13}^{-1}({z_1\over
z_3})R_{23}^{-1}({z_2\over z_3}).}
The general case consists of the $R$-matrices arising on one side in the
normal ordering of the product of the arguments of $\chi$, with factors
$q^{-2c}$, followed by the inverse pattern of $R$-matrices without the
$q^{-2c}$ factor.

%\bibliographystyle{unsrt}
%\bibliography{biblio}

\end{document}